\documentclass[leqno,twoside, 12pt]{amsart}
\usepackage{amssymb, latexsym, esint}
\usepackage{color}

\usepackage{amsthm}
\usepackage{amssymb,amsfonts}
\usepackage[all,arc]{xy}
\usepackage{enumerate}
\usepackage{mathrsfs}
\usepackage{amsmath}
\usepackage{verbatim}

\usepackage[left=2.1cm,right=2.1cm,top=2.9cm,bottom=2.9cm]{geometry}

\theoremstyle{plain}
\numberwithin{equation}{section} \numberwithin{figure}{section}

\newtheorem{theorem}{Theorem}[section]
\newtheorem{lemma}[theorem]{Lemma}
\newtheorem{proposition}[theorem]{Proposition}
\newtheorem{corollary}[theorem]{Corollary}
\newtheorem{definition}[theorem]{Definition}
  \usepackage{epsfig}  
\usepackage{graphicx}  \usepackage{epstopdf}
\theoremstyle{definition}
\newtheorem{remark}[theorem]{Remark}

\numberwithin{equation}{section}

\usepackage[bookmarksnumbered,colorlinks, linkcolor=blue, citecolor=red, pagebackref, bookmarks, breaklinks]{hyperref}

\begin{document}

\title{Embeddings for Orlicz potential spaces with applications to mixed local-nonlocal problems}

\author{Pablo Ochoa}

\address{P. Ochoa. \newline Universidad Nacional de Cuyo, Fac. de Ingenier\'ia. CONICET. Universidad J. A. Maza\\Parque Gral. San Mart\'in 5500\\
Mendoza, Argentina.}
\email{pablo.ochoa@ingenieria.uncuyo.edu.ar}

\author{Ariel Salort}
\address{A. Salort. \newline Departamento de Matemáticas Aplicadas y Ciencia de Datos, 
Universidad San Pablo-CEU, CEU Universities, 
Julián Romea 2, 28003 Madrid, Spain}
\email{ariel.salort@ceu.es}

\parskip 3pt

\maketitle

\begin{abstract}
In this paper, we provide further insight into potential spaces within the nonstandard framework of Orlicz functions. More precisely, we  establish continuous and compact embeddings of potential spaces in $\mathbb{R}^n$ into both unweighted and weighted Orlicz spaces. Special attention is devoted to the radial setting, for which compactness results are obtained through Strauss- and Ni-type inequalities. As an application, we prove the existence of positive radial solutions to a class of Hénon-type problems involving both local and nonlocal contributions in the Orlicz framework. 
\end{abstract}

\section{Introduction}

The \emph{fractional Orlicz-Sobolev spaces} $W^{s,A}(\mathbb R^n)$ introduced in \cite{ACPS, FBS} are the natural generalizations of the classical fractional order Sobolev spaces $W^{s,p}(\mathbb R^n)$, $p>1$.  Given an N-function $A$ and a fractional parameter $s\in (0,1)$, these spaces are defined as
$$
W^{s,A}(\mathbb R^n)=\left\{u\colon \mathbb R^n\to \mathbb R: \int_{\mathbb R^n} A(|u|)\,dx + \iint_{\mathbb R^n\times \mathbb R^n} A\left(\frac{|u(x)-u(y)|}{|x-y|^{s}} \right)\frac{dxdy}{|x-y|^n} <\infty \right\}.
$$
For properties and further research on these spaces we refer to \cite{ACPS, ACPS1,ABS, BOT, FBS}.

In \cite{OS}, a characterization of the fractional Orlicz-Sobolev spaces $W^{s,A}(\mathbb R^n)$ via Bessel potentials is introduced. More precisely, the \emph{Bessel-Orlicz potential space} $H^{s, A}(\mathbb R^n)$ is defined as the set
$$
H^{s, A}(\mathbb{R}^n):=\left\lbrace u\in L^A(\mathbb{R}^n): u= G_s*f, \text{ for some }f\in L^A(\mathbb{R}^n)\right\rbrace,
$$
endowed with the norm $\|u\|_{H^{s, A}(\mathbb{R}^n)}:=\|f\|_{L^A(\mathbb{R}^n)}$, where $\|\cdot\|_{L^A(\mathbb{R}^n)}$ denotes the Luxemburg norm associated with the Orlicz space $L^A(\mathbb{R}^n)$. These spaces are a direct generalization of the spaces studied in \cite{AS, ASmith, Triebel}. In particular, as proved in \cite[Theorems 3.7 and 3.9]{OS}, for $0 < s < s' < 1$ the following continuous embeddings hold:
\begin{equation}\label{cont embdd}
H^{s',A}(\mathbb R^n) \hookrightarrow W^{s,A}(\mathbb R^n) \quad \text{and} \quad W^{s',A}(\mathbb R^n) \hookrightarrow H^{s,A}(\mathbb R^n).
\end{equation}
Furthermore, both spaces coincide in the limit case $s = 1$, i.e., $H^{1,A}(\mathbb R^n) = W^{1,A}(\mathbb R^n)$ with equivalent norms.

Regarding embeddings in Orlicz spaces, in \cite[Theorem 6.1]{ACPS} the authors established the optimal Orlicz target space for the Orlicz-Sobolev embedding. They prove that the continuous embedding $W^{s,A}(\mathbb{R}^n) \hookrightarrow L^{A_{n/s}}(\mathbb{R}^n)$ holds, where  $A_{n/s}$ is defined as (see Section \ref{N functions} for further details)
\begin{equation*}
A_{n/s}(t)=A(H^{-1}(t)), \quad \, H(t)=\left(\int_0^t \left( \frac{\tau}{A(\tau)} \right)^\frac{s}{n-s} \,d\tau \right)^\frac{n-s}{n}, \qquad t\geq 0.
\end{equation*} 
Using the aforementioned embedding, and  combining the potential characterization of Orlicz spaces from \cite{OS} with the boundedness of the Riesz potential established in \cite{BCS}, we obtain in Theorem \ref{teo1} the optimal target space for the embedding of the potential space $H^{s,A}(\mathbb R^n)$ into $L^{A_{n/s}}(\mathbb R^n)$. This result is  applied in Theorem \ref{bounded domains} to get embeddings of $H^{s,A}_0(\Omega)$ into an unweighted Orlicz space $L^B(\Omega)$ for a suitable N-function growing slower than $A$, being $\Omega$ a bounded domain.

In the study of nonlinear elliptic equations on $\mathbb{R}^n$, a fundamental mathematical obstacle is the lack of compactness in standard Sobolev embeddings. 
However, restricting the functional to radial functions introduces a natural spatial decay that recovers compactness. In this direction, Strauss in 
\cite{Strauss} establishes a pointwise decay of radial functions in $H^1(\mathbb R^n)$ bounded away from the origin. Building upon this principle, in 
\cite{Ni} Ni obtained sharp pointwise bounds which became critical for establishing the asymptotic behavior and uniqueness of positive solutions. These type of estimates were generalized to the potential spaces settings, see for instance \cite{EP, Triebel}. More precisely, if $u\in H_{rad}^{s,p}(\mathbb R^n)$ with $n>1$ and $1<p<\infty$ then
\begin{align*}
|u(x)| \leq 
\begin{cases}
c_n |x|^{\frac{1-n}{p}}\|u\|_{H^{s,p}(\mathbb R^n)} &\quad 1/p < s < n \quad\text{(Strauss' estimate)}\\
C_n |x|^{s-\frac{n}{p}}\|u\|_{H^{s,p}(\mathbb R^n)} &\quad 1/p < s < n/p \quad \text{(Ni's estimate)}.
\end{cases}
\end{align*}
In the context of potential spaces in the Orlicz settings, assuming that $A$ satisfies the condition $1<p^-\leq \frac{t A'(t)}{A(t)}\leq p^+<\infty$ for some constants $p^\pm$, in \cite[Theorem 4.3]{OS} we obtained the following Strauss type inequality for $u\in H^{s,A}_{rad}(\mathbb R^n)$:

(\emph{Strauss' estimate})
If $\frac{1}{p^-} < s\leq 1$, there exists $C>0$ such that
$$
|u(x)| \le C \frac{|x|^{1-n}}{\hat{A}^{-1}(|x|^{1-n})} \|u\|_{H^{s,A}(\mathbb{R}^n)}.
$$
To complement this analysis, in Theorem \ref{Ni Lemma} we establish a Ni-type decay inequality for radial functions in the potential space $H^{s,A}(\mathbb R^n)$:

(\emph{Ni's estimate})
If  $\frac{1}{p^-} < s<\frac{n}{p^+}$, there exists $C>0$ such that
$$
|u(x)|\leq C |x|^s A^{-1}(|x|^{-n}) \|u\|_{H^{s,A}(\mathbb R^n)}.
$$

By using these inequalities, we establish in Theorem \ref{teo.rad.compacto} that the embedding 
$$H_{rad}^{s,A}(\mathbb R^n) \subset  L^B(\mathbb R^n, |x|^\alpha\,dx)$$ 
is compact for appropriate relations between the $\Delta_2$-exponents  of the Orlicz functions:
  \begin{equation}\label{cond powers}
  (p_B^--p_A^+)\left(\frac{1-n}{p_{A}^+}\right)+\alpha<0, \qquad \text{and} \qquad p_{B}^+< \frac{(\alpha+n)p_{A}^-}{n-sp_{A}^-}.
  \end{equation}
 These conditions are extensions of the standard-growth case (see for instance \cite{Sin} and \cite[Theorem 3.7]{OS2}). Consequently, this allow us to derive  in Theorem  \ref{teo.sob.compacto} that the embedding 
$$
W^{s,A}(\mathbb R^n) \subset L^B(\mathbb R^n, |x|^\alpha\,dx)
$$
is also compact.

As an application of the above results, we consider a  general class of mixed local and nonlocal Hénon type problems of the form:
\begin{equation} \label{eq general}
\begin{cases}
\mathcal{L}_{s,A,\gamma}\, u + \dfrac{a(|u|)}{|u|}u =|x|^\alpha \dfrac{b(|u|)}{|u|}u  \quad \text{in }\mathbb{R}^n, \\
u >0,
\end{cases}
\end{equation}where $u$ belongs to a suitable functional space of radial functions $\mathcal{W}^{s,A}(\mathbb{R}^n)$  (see \eqref{space}). Given $\gamma\in [0, 1]$, $s\in (0, 1]$ (depending on $\gamma$), the operator $\mathcal{L}_{s,A,\gamma}$ is 
\begin{equation*}
\mathcal{L}_{s,A,\gamma}\, u := \begin{cases}
(-\Delta_a) u, \quad & \text{if}\quad \gamma=1,\\
\gamma (-\Delta_a) u+ (1-\gamma)(-\Delta_a)^s u, \quad & \text{if}\quad \gamma\in (0,1), \, s\in (0,1),\\(-\Delta_a)^s u, \quad & \text{if}\quad \gamma=0,  \, s\in (0,1),
\end{cases}
\end{equation*}where $(-\Delta_a)^s$ is the fractional $a-$Laplacian of order $s\in(0,1]$ (see Section \ref{application} for details). The classical H\'enon's equation \cite{H} with Dirichlet boundary conditions is 
\begin{equation}\label{H eq}
-\Delta u= |x|^\alpha u^{q-1}  \text{ in }B, \qquad 
u>0  \text{ in }B, \qquad 
u=0  \text{ on }\partial B,
\end{equation}where $B$ is the unit ball in $\mathbb{R}^n$, $q>1$ and $\alpha>0$. This  was used to investigate numerically the stability of spherical steady-state stellar systems with respect to spherical disturbances, and it has been the starting point of extensive research. A remarkable progress in the study of \eqref{H eq} was made in the work \cite{Ni}, where it was proved that a  radial solution exists for \eqref{H eq} if
\begin{equation}\label{Ni cond}
2<q<2_\alpha^* := 2^{*}+\dfrac{2\alpha}{n-2} = \frac{2(\alpha+n)}{n-2}. 
\end{equation} Other extensions of problem \eqref{H eq} have been studied in recent years, including nonlinear operators \cite{Nazarov, DP}, different boundary conditions \cite{GS}, nonlocal operators \cite{SW, BQ}, and mixed local–nonlocal operators \cite{SV, OS2}, to mention just a few. 

Further generalizations  to the fractional Laplacian have been obtained in the unbounded setting, where the ball is replaced by the whole space $\mathbb{R}^n$. In \cite{BQ} it was shown that the equation
$$
(-\Delta)^s u = |x|^\alpha u^{q-1} \text{ in } \mathbb{R}^n, \qquad u>0, 
$$
admits no positive solutions in the range $~2<q<2^*_{s,\alpha}:=\frac{2(n+\alpha)}{n-2s}$, where $\alpha>-2s$, due to the loss of mass at infinity. In order to obtain existence of solutions in   $\mathbb{R}^n$, a confining potential can be added to prevent the loss of mass at infinity (see for instance \cite{Sin}, \cite{Ma} or \cite{OS2}). Following these references, in Section \ref{application} we apply the continuous and compact embeddings developed in this work  to employ  critical point theory to extend the above existence results for the  H\'enon-type model \eqref{eq general}. As a result, in Theorem \ref{existence theorem} we obtain  that under condition \eqref{cond powers}, problem \eqref{eq general} has a radial solution. Observe that \eqref{cond powers} constitutes a generalization of the Ni's conditions \eqref{Ni cond}.

The article is organized as follows. In Section \ref{sec.prelim} we introduce the basic notation and known results on fractional Orlicz-Soblev spaces, rearrangements and the Bessel kernel. Section \ref{sec.4} are devoted to prove our continuous embedding results while that Section \ref{sec.5} contains our compact embedding results for potential spaces and radial functions. Finally, in Section \ref{application} we address an application related with the H\'enon problem for a mixed local-nonlocal operator in $\mathbb R^n$.

\section{Preliminaries} \label{sec.prelim}

In this section, we introduce basic definitions and preliminary results  used in the paper.

\subsection{$N$-functions}\label{N functions}
\begin{definition}
We say that $A:[0, \infty) \to [0, \infty)$ is an $N$-function if satisfies the following conditions
\begin{enumerate}
\item[(i)] $A$ is convex, increasing, and $A(0) = 0$.
\item[(ii)] $A$ is super-linear at zero and at infinity, that is 
$$
\lim_{t\rightarrow 0+} \frac{A(t)}{t}=0\quad \text{and}\quad \lim_{t\rightarrow \infty} \frac{A(t)}{t}=\infty.
$$
\end{enumerate}
\end{definition}
Moreover, an $N$-function $A$ admits the representation
$$
A(t)=\int_0^t a(\tau)\,d\tau,
$$
where the function $a$ is right-continuous for $t \geq 0$, positive for $t>0$, non-decreasing and
satisfies the conditions 
$$
a(0)=0, \qquad \lim_{t\to\infty} a(t)=\infty.
$$
We extend $A$ to $\mathbb R$ as an even function, and we assume without loss of generality that $A(1)=1$.

\noindent  Associated to $A$ is  the complementary N-function function   to it which is defined as follows:
\begin{equation}\label{Gcomp}
\hat{A} (t) := \sup \left\lbrace tw-A(w) \colon w>0 \right\rbrace .
\end{equation}
The following Young-type inequality holds
\begin{equation}\label{2.5}
t_1t_2 \leq A(t_1)+\hat{A} (t_2) \text{  for every } t_1,t_2 \geq 0.
\end{equation}

An $N$-function $A$ satisfies the $\Delta_2$ condition if there exists $C>2$ such that 
$$
A(2t)\leq CA(t) \quad \text{ for all } t\geq 0.
$$

\noindent By \cite[Theorems 4.1 and 4.3, Chapter 1]{KR},  an $N$-function $A$ and its complementary N-function $\hat A$ satisfy the $\Delta_{2}$ condition if and only if there is $1<p^-\leq p^+ <\infty$ such that
\begin{equation}\label{eq.p}
p^-\leq \frac{ta(t)}{A(t)} \leq p^+, ~~~~~\forall\, t>0.
\end{equation}
In particular, \eqref{eq.p} implies that
\begin{align} \label{potencias}
\min\{t^{p^-}, t^{p^+}\}&\leq A(t) \leq \max\{t^{p^-}, t^{p^+}\}, \quad t\geq 0.
\end{align}

For further reference, we next define a special $N$-function. Let $A$ be an $N$-function such that \begin{equation} \label{cond.A}
\int^\infty \left( \frac{t}{A(t)}\right)^\frac{s}{n-s}\,dt =\infty, \qquad  \int_0 \left( \frac{t}{A(t)}\right)^\frac{s}{n-s}\,dt <\infty,
\end{equation}
then the $N$-function $A_{n/s}$ is defined as (see \cite{ACPS})
\begin{equation} \label{A.ns}
A_{n/s}(t)=A(H^{-1}(t)), \qquad t\geq 0.
\end{equation}
Here, the auxiliary function $H$ is given by
$$
H(t)=\left(\int_0^t \left( \frac{\tau}{A(\tau)} \right)^\frac{s}{n-s} \,d\tau \right)^\frac{n-s}{n}, \qquad t\geq 0.
$$

As we will see in the next sections, $A_{n/s}$ is the optimal critical $N$-function in the embedding of Orlicz-Sobolev spaces into Orlicz spaces.

\subsection{Orlicz-Sobolev spaces} Given an N-function $A$, we define the \emph{Orlicz space}
\begin{align*}
L^A(\mathbb R^n) = \{ u\colon  \mathbb R^n \to \mathbb R : \rho_A(u)<\infty\},
\end{align*}
and the \emph{fractional Orlicz-Sobolev space} of order $s\in (0,1]$ 
$$
W^{s,A}(\mathbb R^n)=\{ u \colon \mathbb R^n \to \mathbb R : \rho_A(u)+\rho_{s,A}(u) < \infty\},
$$
where $\rho_A$ and $\rho_{s,A}$ are the modulars given by
\begin{align*}
\rho_A(u):=\int_{\mathbb R^n} A(|u|)\,dx, \qquad 
\rho_{s,A}(u):=
\begin{cases}
\displaystyle\iint_{\mathbb R^n\times\mathbb R^n} A(|D^s u|)\,d\mu, & \text{ if } s\in(0,1),\\
\displaystyle\int_{\mathbb R^n} A(|\nabla u|)\,dx & \text{ if } s=1,
\end{cases}
\end{align*}
where $d\mu=|x-y|^{-n}dxdy$ and $D^s u=\frac{u(x)-u(y)}{|x-y|^s}$.

The spaces $L^A(\mathbb R^n)$ and $W^{s,A}(\mathbb R^n)$ are endowed with the Luxemburg norms defined as
$$
\|u\|_{L^A(\mathbb R^n)} = \inf\left\{\lambda>0 \colon \rho_A\left(\frac{u}{\lambda} \right)\leq 1 \right\}
$$
and
$$
\|u\|_{W^{s,A}(\mathbb R^n)} = \|u\|_{L^A(\mathbb R^n)} + [u]_{W^{s,A}(\mathbb R^n)},
$$
being $[u]_{W^{s,A}(\mathbb R^n)}$ the so-called Gagliardo-like seminorm defined as
$$
[u]_{W^{s,A}(\mathbb R^n)} = \inf\left\{ \lambda >0  \colon \rho_{s,A}\left(\frac{u}{\lambda} \right)\leq 1\right\}.
$$ 
 The following H\"older's type inequality holds for Orlicz spaces (see \cite[Theorem 3.7.5]{kufner}).
\begin{lemma}
Let $A$ be an N-function with complementary function $\hat A$. If $u\in L^A(\mathbb R^n)$ and $v\in L^{\hat A}(\mathbb R^n)$ then $uv\in L^1(\mathbb R^n)$ and
$$
\int_{\mathbb R^n} |uv|\,dx \leq 2 \|u\|_{L^A(\mathbb R^n)} \|v\|_{L^{\hat A}(\mathbb R^n)}.
$$
\end{lemma}

 In \cite[Theorem 4.1]{BC} it is proved the following inequality for convolutions in Orlicz spaces.
 \begin{lemma} \label{lema.1}
 Let $A$ be an N-function. Then
 $$
 \|u*v\|_{L^A(\mathbb R^n)} \leq \|v\|_{L^1(\mathbb R^n)}\|u\|_{L^A(\mathbb R^n)}
 $$
 for every $u\in L^A(\mathbb R^n)$ and $v\in L^1(\mathbb R^n)$. Moreover,
 $$
 \int_{\mathbb R^n}A(|u*v|)\,dx \leq\int_{\mathbb R^n} A\left(\|v\|_{L^1(\mathbb R^n)} |u|\right)\,dx.
 $$
 \end{lemma}

 \subsection{Bessel kernels} \label{sec.potential.spaces}

Given $s\in \mathbb R$, we consider the Bessel kernel
$$
G_s(x)=\mathcal{F}^{-1}\left((1+|x|^2)^{-s/2}\right),
$$
where $\mathcal{F}(f)$ denotes the Fourier Transform of $f$.  Here, we allow $s\in \mathbb R$ since $\hat G_s(\xi)=(1+|\xi|^2)^{-s/2}$ identified with a tempered distribution for all $s$.   According to \cite[Proposition 2]{S} and \cite[pages 13 and 27]{ASmith}, we recall some useful properties of $G_s$.
\begin{lemma} \label{lema.G}
The Bessel kernel satisfies the following properties:
\begin{itemize}
\item[i)] $G_s>0$, $G_s\in L^1(\mathbb{R}^n)$ for all $s>0$ and $\|G_s\|_{L^1(\mathbb R^n)}=1$,
\item[ii)] $\nabla G_s\in L^1(\mathbb{R}^n)$ for all $s>1$,
\item[iii)] $G_s(x) = G_s(-x)$ for all $x\in \mathbb R^n$,
\item[iv)] $G_s(x)$ is radial and radially decreasing.
\end{itemize}
\end{lemma}

\subsection{Rearrangements}

We recall some definitions concerning  rearrangement for a measurable function $u:\mathbb R^n \to \mathbb{R}$. The distribution function $\mu$ of $u$ is defined as
$$\mu(t)=|\{x\in \mathbb R^n: |u(x)|> t\}|,\quad \text{for }t>0.$$
We observe that $\mu$ is non-increasing and right-continuous. 

The decreasing rearrangement of $u$, denoted by $u^{*}$,  is the function
$$u^{*}(s)= \inf\{\lambda\geq 0: \mu(\lambda)\leq s\}, \quad \text{for }s\geq 0.$$The symmetric rearrangement of $u$, denoted by $u^{\bigstar}$, is the function
$$u^{\bigstar}(x)=u^{*}(\omega_n|x|^{n}), \quad x\in \mathbb R^n,$$where $\omega_n$ is the volume of the $n$-dimensional unit ball. The following integral representations are valid for $u^{*}$ and $u^{\bigstar}$:
$$
u^{*}(s)=\int_{0}^{\infty}\chi_{[0, \mu(t))}(s)\,dt,
$$
$$
u^{\bigstar}(x)=\int_0^{\infty}\chi_{\{y:\, \omega_n |y|^n < \mu(t)\}}(x)\,dt=\int_0^{\infty}\chi_{\{y:\, \omega_n |y|^n < \mu(t)\}^{*}}(x)\,dt,
$$
where $\{y:\, \omega_n |y|^n < \mu(t)\}^{*}$ is the symmetric rearrangement of the set $\{y:\, \omega_n |y|^n < \mu(t)\}$, that is, it is the open ball centered at the origin whose volume is that of $\{y:\, \omega_n |y|^n < \mu(t)\}.$
 
\section{Continuous embeddings of $H^{s, A}$ into Orlicz spces } \label{sec.4}

By combining the potential characterization of Orlicz spaces from \cite{OS} with the boundedness of the Riesz potential established in \cite{BCS}, we obtain a very simple proof of the optimal target space for the embedding of the potential space $H^{s,A}$.

\begin{theorem} \label{teo1}
Let $A$ be an $N$-function satisfying \eqref{cond.A}, and let $A_{n/s}$ be defined as in \eqref{A.ns}. Then
$$
H^{s,A}(\mathbb R^n) \hookrightarrow L^{A_{n/s}}(\mathbb R^n)
$$
and there exists $C>0$ such that
$$
\|u\|_{L^{A_{n/s}}(\mathbb R^n)} \leq C\|u\|_{H^{s,A}(\mathbb R^n)}.
$$
\end{theorem}

\begin{proof}

Let $u \in H^{s,A}(\mathbb R^n)$. By definition, there exists a function $f \in L^A(\mathbb R^n)$ such that $u = G_s * f$ and $\|u\|_{H^{s,A}(\mathbb R^n)} = \|f\|_{L^A(\mathbb R^n)}$. Since the Bessel kernel $G_s$ and the Riesz kernel $I_s$ satisfy the pointwise inequality $G_s(x) \leq I_s(x)$ for all $x \in \mathbb R^n \setminus \{0\}$, we obtain
$$
|(G_s * f)(x)| \leq (G_s * |f|)(x) \leq (I_s * |f|)(x) \quad \text{for a.e. } x \in \mathbb R^n.
$$
By the monotonicity of the Orlicz norm, this pointwise inequality implies
$$
\|G_s * f(x)\|_{L^{A_{n/s}}(\mathbb R^n)} \leq\| I_s * |f|\|_{L^{A_{n/s}}(\mathbb R^n)}.
$$
According to \cite[Theorem 3.6]{BCS}, the Riesz potential operator is bounded from $L^A(\mathbb R^n)$ to $L^{A_{n/s}}(\mathbb R^n)$. Thus, there exists a constant $c = c(s, n) > 0$ such that
$$
\|I_s * |f|\|_{L^{A_{n/s}}(\mathbb R^n)} \leq c \||f|\|_{L^A(\mathbb R^n)} = c \|f\|_{L^A(\mathbb R^n)}.
$$
Combining these inequalities, we conclude that
$$
\|u\|_{L^{A_{n/s}}(\mathbb R^n)} =  \|G_s * f(x)\|_{L^{A_{n/s}}(\mathbb R^n)} \leq c\|f\|_{L^A(\mathbb R^n)} = c\|u\|_{H^{s,A}(\mathbb R^n)}.
$$
which completes the proof.
\end{proof}

\begin{corollary}\label{continuity}
Let $B$ be an N-function. Assume there exist constants $c, t_0 > 0$ such that:
\begin{itemize}
\item[i)] $B(t)\leq A(ct) \quad \text{ for all }0\leq t\leq t_0$,
\item[ii)] $B(t)\leq A_{n/s}(ct) \quad \text{ for all } t\geq t_0$.
\end{itemize}
Then, the embedding $H^{s,A}(\mathbb R^n) \hookrightarrow L^B(\mathbb R^n)$ is continuous. In particular, there exists a constant $C > 0$ such that
$$
\|u\|_{L^B(\mathbb R^n)} \leq C \|u\|_{H^{s,A}(\mathbb R^n)} \quad \text{for all } u \in H^{s,A}(\mathbb R^n).
$$
\end{corollary}

\begin{proof}
Since $u \in H^{s,A}(\mathbb R^n)$, there exists $f \in L^A(\mathbb R^n)$ such that $u = G_s * f$ and $\|u\|_{H^{s,A}(\mathbb R^n)} = \|f\|_{L^A(\mathbb R^n)}$. First, since N-functions are nonnegative, conditions \textup{(i)} and \textup{(ii)} imply that 
$$
B(t) \leq A(ct) + A_{n/s}(ct), \qquad t\geq 0.
$$
By Theorem \ref{teo1}, there exists a constant $C_1 > 0$ such that
\begin{equation} \label{cota.0}
\|u\|_{L^{A_{n/s}}(\mathbb R^n)} \leq C_1 \|u\|_{H^{s,A}(\mathbb R^n)}.
\end{equation}
Moreover, since $G_s \in L^1(\mathbb R^n)$ by Lemma \ref{lema.G}, Young's inequality for convolution in Orlicz spaces yields
\begin{equation} \label{cota.1}
\|u\|_{L^A(\mathbb R^n)} = \|G_s * f\|_{L^A(\mathbb R^n)} \leq \|G_s\|_{L^1(\mathbb R^n)} \|f\|_{L^A(\mathbb R^n)} =  \|u\|_{H^{s,A}(\mathbb R^n)}.
\end{equation}
Next, let $\lambda = \max\left\{ \|u\|_{L^A(\mathbb R^n)}, \|u\|_{L^{A_{n/s}}(\mathbb R^n)} \right\}$. By the definition of the Luxemburg norm, we have
\begin{equation} \label{relac.1}
\int_{\mathbb R^n} A\left(\frac{|u(x)|}{\lambda}\right) \,dx \leq 1 \quad \text{and} \quad \int_{\mathbb R^n} A_{n/s}\left(\frac{|u(x)|}{\lambda}\right) \,dx \leq 1.
\end{equation}
Using the convexity of the N-functions together with \eqref{relac.1}, we obtain
\begin{align*}
\int_{\mathbb R^n} B\left(\frac{|u(x)|}{2\lambda c} \right)\,dx &\leq 
\int_{\mathbb R^n} A\left(\frac{c|u(x)|}{2\lambda c} \right)\,dx +
\int_{\mathbb R^n} A_{n/s}\left(\frac{c|u(x)|}{2\lambda c} \right)\,dx\\
&\leq
\frac12 \int_{\mathbb R^n} A\left(\frac{|u(x)|}{\lambda } \right)\,dx +
\frac12 \int_{\mathbb R^n} A_{n/s}\left(\frac{|u(x)|}{\lambda } \right)\,dx=1.
\end{align*}
Therefore, by the definition of the Luxemburg norm combined with \eqref{cota.0} and \eqref{cota.1} gives
\[
\|u\|_{L^B(\mathbb R^n)} \leq 2c \max\left\{ \|u\|_{L^A(\mathbb R^n)}, \|u\|_{L^{A_{n/s}}(\mathbb R^n)} \right\} \leq 2c C_1 \|u\|_{H^{s,A}(\mathbb R^n)},
\]
which completes the proof.
\end{proof}

Using \cite[Theorem 3.9]{OS}, we easily obtain the following result.
\begin{corollary}
Let $s'>s>0$. Then $W^{s',A}(\mathbb R^n) \hookrightarrow L^{A_{n/s}}(\mathbb R^n)$, and there exists $C>0$ such that
$$
\|u\|_{L^{A_{n/s}}(\mathbb R^n)} \leq C \|u\|_{W^{s',A}(\mathbb R^n)}
$$

\end{corollary}

In the end of the section, we consider embedding of potential spaces into Orlicz spaces in bounded domains. 

Let $\Omega\subset \mathbb R^n$ be an open bounded set. Define
$$
H^{s,A}_0(\Omega) := \{u \in H^{s,A}(\mathbb R^n)  \colon u=f*G_s \text{ with } f\in L^A(\Omega) \text{ and } f=0 \text{ in } \mathbb R^n\setminus \Omega\}.
$$
\begin{theorem}\label{bounded domains}
If there are $c>0$ and $t_0>0$ such that $B(t)\leq A(ct)$ for all $t\geq t_0$, then $H_0^{s,A}(\Omega)\subset L^B(\Omega)$.
\end{theorem}

To prove Theorem \ref{bounded domains}, we will state a preliminary inequality for rearrangements.  Given a function $u$, we recall that  $u^*$ denotes its decreasing rearrangment, and $u^\bigstar$ is the symmetric rearrangement of $u$.

We quote the following rearrangement inequality (see for instance \cite[Theorem 3.7]{LL}).

\begin{proposition}[Riesz rearrangement inequality] \label{riesz}
Any three non-negative functions $f : \mathbb{R}^n \to \mathbb{R}^+$, $g : \mathbb{R}^n \to \mathbb{R}^+$ and $h : \mathbb{R}^n \to \mathbb{R}^+$ satisfy the inequality
$$
\iint_{\mathbb{R}^n \times \mathbb{R}^n} f(y)g(x - y)h(x) \, dy \, dx \le \iint_{\mathbb{R}^n \times \mathbb{R}^n} f^{\bigstar}(y)g^{\bigstar}(x - y)h^{\bigstar}(x) \, dy \, dx,
$$
where $f^{\bigstar} : \mathbb{R}^n \to \mathbb{R}^+$, $g^{\bigstar} : \mathbb{R}^n \to \mathbb{R}^+$ and $h^{\bigstar} : \mathbb{R}^n \to \mathbb{R}^+$ are the symmetric decreasing rearrangements of the functions $f$, $g$ and $h$ respectively.
\end{proposition}

\begin{theorem} \label{thm.PS}
Let $u\in H^{s,A}(\mathbb R^n)$ be such that $u=f*G_s$ for some $f\in L^A(\mathbb R^n)$. Define $v=f^{\bigstar}*G_s$. Then
$$
\int_0^t u^*(r) \,dr \leq \int_0^t v^*(r)\,dr.
$$
\end{theorem}

\begin{proof}
Let $u\in H^{s,A}(\mathbb R^n)$. Then there is $f\in L^A(\mathbb R^n)$ such that $u=f*G_s$.

Since $G_s\ge0$, $|u|\le|f|*G_s$ pointwise, and the decreasing rearrangement is monotone, so $u^*\le(|f|*G_s)^*$. Hence it suffices to prove the theorem with $u$ replaced by $\tilde u:=|f|*G_s\ge0$, whose associated function is the same $v=f^{\bigstar}*G_s$. From now on we assume $f\ge0$ and $u=f*G_s\ge0$.

We consider the function $f^{\bigstar}$, the symmetric rearrangement of $f$. Then, $\|f^{\bigstar}\|_{L^A(\mathbb R^n)}=\|f\|_{L^A(\mathbb R^n)}$.

Since $f^{\bigstar}\in L^A(\mathbb R^n)$, we can define the function $v=f^{\bigstar} * G_s \in H^{s,A}(\mathbb R^n)$, and
$$
\|v\|_{H^{s,A}(\mathbb R^n)} = \|f^{\bigstar}\|_{L^A(\mathbb R^n)}=\|f\|_{L^A(\mathbb R^n)} = \|u\|_{H^{s,A}(\mathbb R^n)}
$$
Observe that, since both $G_s$ is positive, radial and decreasing, $v$ is radial and strictly decreasing. From Proposition \ref{riesz}, and since $G_s^{\bigstar} =G_s$,  we get
$$
\iint_{\mathbb{R}^n \times \mathbb{R}^n} f(y)G_s(x - y)h(x) \, dy \, dx \le \iint_{\mathbb{R}^n \times \mathbb{R}^n} f^{\bigstar}(y)G_s(x - y)h^{\bigstar}(x) \, dy \, dx.
$$
Let $E\subset \mathbb R^n$ be an open set such that $|E|=t$, and take $h(x)=\chi_E(x)$. In this case we have that $h^{\bigstar}(x)=\chi_{B_{R_t}}(x)$, where $R_t$ is such that $|B_{R_t}|=t$. The previous expression gives that
$$
\int_E \int_{\mathbb{R}^n} f(y)G_s(x - y) \, dy \, dx \le  \int_{B_{R_t}} \int_{\mathbb{R}^n} f^{\bigstar}(y)G_s(x - y)\, dy \, dx,
$$
which can be written as
$$
\int_E f*G_s(x) \, dx \le  \int_{B_{R_t}} f^{\bigstar}*G_s(x) \, dx,
$$
that is, 
\begin{equation} \label{des.E}
\int_E u(x) \, dx \le  \int_{B_{R_t}} v(x)  \, dx.
\end{equation}
By \cite[Lemma 2.5]{BS},
$$
\int_0^t u^*(r) \, dr = \sup_{|E|=t} \int_E u(x) \, dx \, \text{ and }\int_0^t v^*(r) \, dr = \sup_{|E|=t} \int_E v(x) \, dx.
$$
Taking supremum over all $E$ such that $|E|=t$ in \eqref{des.E}, we obtain that
$$
\int_0^t u^*(r) \, dr \le  \int_{B_{R_t}} v(x) 
 \, dx \leq  \sup_{|E|=t} \int_E v(x) \,dx = \int_0^t v^*(r) \, dr.
 $$
This completes the proof.
\end{proof}

 Now, we are ready to prove Theorem \ref{bounded domains}.

\begin{proof}[Proof of Theorem \ref{bounded domains}]
Let $u\in H^{s,A}_0(\Omega)$. Then there exists $f\in L^A(\mathbb R^n)$ such that $f=0$ in $\mathbb R^n \setminus \Omega$. Then $f^{\bigstar} \in L^A(\mathbb R^n)$ and $f^{\bigstar}=0$ in $\mathbb R^n \setminus \Omega^*$ for some ball $\Omega^*$. By Theorem \ref{thm.PS}, 
$$
\int_0^t u^*(r) \, dr \le \int_0^t v^*(r)\,dr.
$$
Therefore, by \cite[Corollary 4.7]{BS}, since $L^B(\mathbb R^n)$ is a rearrangement-invariant space, if $v\in L^B(\Omega)$ then $u\in L^B(\Omega)$ and 
\begin{align}
\|u\|_{L^B(\Omega)} \leq \|u\|_{L^B(\mathbb R^n)} \leq \|v\|_{L^B(\mathbb R^n)} \leq  C \|f\|_{L^B(\mathbb R^n)} \leq C(|\Omega|) \|f\|_{L^A(\mathbb R^n)} =  C \|u\|_{H^{s,A}(\mathbb R^n)}.
\end{align}
where we have used Young's to bound
\begin{align*}
\|v\|_{L^B(\Omega)}&\leq \|v\|_{L^B(\mathbb R^n)} \leq  \|f^{\bigstar}*G_s\|_{L^B(\mathbb R^n)}\leq \|f^{\bigstar}\|_{L^B(\mathbb R^n)} \|G_s\|_{L^1(\mathbb R^n)}\\
&= 
\|f^{\bigstar}\|_{L^B(\Omega^*)} \|G_s\|_{L^1(\mathbb R^n)}
= 
\|f\|_{L^B(\mathbb R^n)} \|G_s\|_{L^1(\mathbb R^n)},
\end{align*}and the embedding $L^A(\mathbb R^n)\subset L^B(\mathbb R^n)$. This completes the proof.
\end{proof}

\section{Compact embeddings for  potential spaces of radial functions} \label{sec.5}

In this section, we will provide conditions for the compactness of the inclusion of $H^{s, A}$ into Orlicz spaces. First, we will discuss the unweighted case and then we will give conditions for the validity of the compact inclusion into weighted Orlicz spaces.

Consider the following class of radial potential  functions
$$H_{rad}^{s, A}(\mathbb{R}^n)=\left\lbrace u\in H^{s, A}(\mathbb{R}^n): \text{ $u$ is radially symmetric} \right\rbrace.$$

In order to obtain compact embeddings of potential spaces into weighted Orlicz spaces, we will apply results regarding the decay of radial functions. Firstly, in \cite{OS}, we obtained the following Strauss type inequality
\begin{theorem}\label{Strauss Lemma}
Let $A$ be a N-function so that $A$ and $\hat{A}$ satisfy the $\Delta_2$-condition. Let $u \in H_{rad}^{s,A}(\mathbb R^n)$ with $\frac{1}{p^-} < s\leq 1$. Then there exists $C > 0$ depending of $n$, $p^+$ and $p^-$ such that
$$
|u(x)| \le C \frac{|x|^{1-n}}{\hat{A}^{-1}(|x|^{1-n})} \|u\|_{H^{s,A}(\mathbb{R}^n)}.
$$
\end{theorem}

The first compactness result is the following.

\begin{theorem}\label{compact emb}
Assume that $A$ and  $\hat A$ are $N$-function satisfying the $\Delta_2$-condition. Let $\frac{1}{p^-}<s<\frac{n}{p^+}$. Let $B$ be an $N$-function such that
$$\lim_{t\to 0^{+}}\dfrac{B(t)}{A(t)}=0 \quad \text{and}\quad B\ll  A_{n/s}.$$
Then, the embedding $H_{rad}^{s,A}(\mathbb R^n) \hookrightarrow L^B(\mathbb R^n)$ is compact.
\end{theorem}

\begin{proof}
Let $\left\lbrace u_k\right\rbrace_{k\in\mathbb N} \subset H_{rad}^{s, A}(\mathbb R^n)$ be a bounded sequence. Then, up to a subsequence $u_k\rightharpoonup u$. Set $v_k:=u_k-u$. Then $v_k\in H^{s,A}_{rad}(\mathbb R^n)$, $v_k\rightharpoonup0$ and $M:=\sup_k\|v_k\|_{H^{s,A}(\mathbb R^n)}<\infty$. We must show $v_k\to0$ in $L^B(\mathbb R^n)$.
 
\smallskip
\emph{Step 1: uniform bounds.} By Young's inequality for convolutions, $\|v_k\|_{L^A(\mathbb R^n)}\le M$, so
$$\int_{\mathbb R^n}A(|v_k|)\,dx\le C_0:=\max\{M^{p^-},M^{p^+}\}.$$
By Theorem \ref{teo1}, $K:=\sup_k\|v_k\|_{L^{A_{n/s}}(\mathbb R^n)}<\infty$, hence $\int_{\mathbb R^n}A_{n/s}(|v_k|/K)\,dx\le 1$.
 
\smallskip
\emph{Step 2: uniform decay.} By  the Strauss-type Theorem \ref{Strauss Lemma} and the inequality $t\le A^{-1}(t)\hat A^{-1}(t)$,
$$|v_k(x)|\le C\,\frac{|x|^{1-n}}{\hat A^{-1}(|x|^{1-n})}\,M\le C M\,A^{-1}(|x|^{1-n}),\qquad x\neq0.$$
Since $n\ge2$, $A^{-1}(|x|^{1-n})\to0$ as $|x|\to\infty$. Thus, for every $\delta>0$ there is $R=R(\delta)>0$ such that
$$|v_k(x)|\le\delta \quad\text{for a.e. } |x|\ge R, \text{ for all } k.$$
 
\smallskip
\emph{Step 3: local compactness.} By \cite[Theorem 3.7]{OS}, $\{v_k\}$ is bounded in $W^{s',A}(\mathbb R^n)$ for $0<s'<s$, hence in $W^{s',A}(B_R)$. By \cite[Theorem 9.1]{ACPS}, it is precompact in $L^A(B_R)$. Since also $v_k\rightharpoonup0$ in $L^A(B_R)$, every subsequence has a further subsequence converging strongly to $0$; hence $v_k\to0$ in $L^A(B_R)$ and, as $A\in\Delta_2$,
$$\int_{B_R}A(|v_k|)\,dx\to0\quad\text{as }k\to\infty.$$
 
\smallskip
\emph{Step 4: conclusion.} Fix $\lambda>0$ and $\varepsilon>0$. Since $A(\lambda t)\le\max\{\lambda^{p^-},\lambda^{p^+}\}A(t)$,
$$
\frac{B(\lambda t)}{A(t)}=\frac{B(\lambda t)}{A(\lambda t)}\frac{A(\lambda t)}{A(t)}\to0 \quad(t\to0^+).
$$
We have the following situations:

(a)\, Choose $\delta\in(0,1)$ with $B(\lambda t)\le\varepsilon A(t)$ for $0\le t\le\delta$.

(b)\, Since $B\ll A_{n/s}$, choose $T>\delta$ with $B(\lambda t)\le\varepsilon A_{n/s}(t/K)$ for $t\ge T$.

(c)\, For $\delta\le t\le T$, $B(\lambda t)\le B(\lambda T)\le C_*A(t)$, with $C_*:=B(\lambda T)/A(\delta)$.

Take $R=R(\delta)$ from Step 2.

On $\{|x|\ge R\}$ we have $|v_k|\le\delta$, so
$$\int_{|x|\ge R}B(\lambda|v_k|)\,dx\le\varepsilon\int_{\mathbb R^n}A(|v_k|)\,dx\le\varepsilon C_0.$$

On $B_R$, by the three cases above,
$$\int_{B_R}B(\lambda|v_k|)\,dx\le(\varepsilon+C_*)\int_{B_R}A(|v_k|)\,dx +\varepsilon\int_{\mathbb R^n}A_{n/s}\!\Big(\frac{|v_k|}{K}\Big)dx \le(\varepsilon+C_*)\int_{B_R}A(|v_k|)\,dx+\varepsilon.$$
By Step 3 the first term tends to $0$, so $\limsup_k\int_{\mathbb R^n}B(\lambda|v_k|)\,dx\le\varepsilon(C_0+1)$. As $\varepsilon$ is arbitrary, $\int_{\mathbb R^n}B(\lambda|v_k|)\,dx\to0$ for every $\lambda>0$. Taking $\lambda=1/\eta$ gives $\|v_k\|_{L^B(\mathbb R^n)}\le\eta$ for $k$ large, for every $\eta>0$. Hence $u_k\to u$ in $L^B(\mathbb R^n)$.
\end{proof}

To complement our analysis, we now establish a Ni-type decay inequality for radial functions in the potential space.

\begin{theorem}\label{Ni Lemma}
Let $s\in (0, 1]$ and let $A$ be an N-function such that $\frac{1}{p^{-}} <s<\frac{n}{p^+}$. Then, for $u\in H^{s,A}_{rad}(\mathbb R^n)$ 
$$
|u(x)|\leq C |x|^s A^{-1}(|x|^{-n}) \|u\|_{H^{s,A}(\mathbb R^n)}.
$$
where $C>0$ depends on the bounds of $A$, $n$ and $s$.
\end{theorem}

 Before proving Theorem \ref{Ni Lemma}, we prove the following useful estimate for norms.
\begin{lemma} \label{lema.normas}
 Let $r>0$, $S_r=\{r/2<|x|<2r\}$. There is $C_0\ge1$, independent of $r$, such that
\begin{align}
\big\||x|^{s-n}\big\|_{L^{\hat A}(\mathbb R^n\setminus B_r)}&\le C_0\,r^{s}A^{-1}(r^{-n})&&\text{if } s<n/p^+,  s\in (0, 1] ,\label{norma.ext}\\
\|\big||x|-r\big|^{s-1}\|_{L^{\hat A}(S_r)}&\le C_0\,r^{n-1+s}A^{-1}(r^{-n})&&\text{if } 1/p^-<s<1.\label{norma.mid}
\end{align}
\end{lemma}
  
\begin{proof}
Since $t\le A^{-1}(t)\hat A^{-1}(t)$, it gives
\begin{equation}\label{A and inverse hat A}
\frac1{r^nA^{-1}(r^{-n})}\le \hat A^{-1}(r^{-n}),\qquad \hat A\big(\hat A^{-1}(r^{-n})\big)=r^{-n}.
\end{equation}

Let us prove \eqref{norma.ext}. Let $\lambda=C_0r^sA^{-1}(r^{-n})$, $C_0\ge1$. With $x=ry$,
\begin{align*}
\int_{\mathbb R^n\setminus B_r}\hat A\Big(\frac{|x|^{s-n}}{\lambda}\Big)dx
&=r^n\!\int_{|y|>1}\hat A\Big(\frac{|y|^{s-n}}{C_0\,r^nA^{-1}(r^{-n})}\Big)dy
\le r^n\!\int_{|y|>1}\hat A\big(C_0^{-1}|y|^{s-n}\hat A^{-1}(r^{-n})\big)dy\\
&\le C_0^{-\frac{p^+}{p^+-1}}\int_{|y|>1}|y|^{(s-n)\frac{p^+}{p^+-1}}\,dy .
\end{align*}
The last integral is equal to $n\omega_n/\big((n-s)\frac{p^+}{p^+-1}-n\big)$, which is finite provided that $(n-s)\frac{p^+}{p^+-1}>n$, i.e.\ $s<n/p^+$. Taking $C_0\ge1$ large enough, we get the result.
 
 \medskip
 
Let us prove \eqref{norma.mid}. Let $\lambda=C_0\,r^{n-1+s}A^{-1}(r^{-n})$. Using polar coordinates, $\rho=r(1+\tau)$ with $\tau\in(-\tfrac12,1)$ and $d\rho=r\,d\tau$,
\begin{align*}
\int_{S_r}\hat A\Big(\frac{\big||x|-r\big|^{s-1}}{\lambda}\Big)dx
&\le n\omega_n(2r)^{n-1}\int_{r/2}^{2r}\hat A\Big(\frac{|r-\rho|^{s-1}}{\lambda}\Big)d\rho
= C\,r^{n}\int_{-1/2}^{1}\hat A\Big(\frac{|\tau|^{s-1}}{C_0\,r^nA^{-1}(r^{-n})}\Big)d\tau\\
&\le C\,r^{n}\int_{-1/2}^{1}\hat A\big(C_0^{-1}|\tau|^{s-1}\hat A^{-1}(r^{-n})\big)d\tau\\
&\leq
 C\big(C_0^{-\frac{p^+}{p^+-1}}+C_0^{-\frac{p^-}{p^--1}}\big)\int_{-1/2}^{1}\big(|\tau|^{-(1-s)\frac{p^+}{p^+-1}}+|\tau|^{-(1-s)\frac{p^-}{p^--1}}\big)d\tau .
\end{align*}

The integral is finite if and only if $(1-s)\frac{p^-}{p^--1}<1$, i.e., $s>1/p^-$. Taking $C_0$ large the whole expression is lower than 1, which proves \eqref{norma.mid}. 
\end{proof}
 
\begin{proof}[Proof of Theorem \ref{Ni Lemma}]
 
Assume first that $s\in (0, 1)$. Let $0<s'<s$ such that $s'p^{-}>1$. Let $u\in H^{s,A}_{rad}(\mathbb R^n)$. Then  by \cite[Theorem 3.7]{OS} the embedding $H^{s,A}(\mathbb R^n) \subset W^{s',A}(\mathbb R^n)$ holds and so $u\in W^{s',A}_{rad}(\mathbb R^n)$. Since $s'p^->1$, by \cite[Theorem 4.2]{ACPS3} $u$ is almost everywhere equal to a continuous function in $\mathbb R^n\setminus \{0\}$.

Since  $u\in H^{s,A}_{rad}(\mathbb R^n)$, by \cite[Lemma 4.1]{OS} there is $f\in L^A(\mathbb R^n)$ radial such that
 
\begin{align*}
u(x)= \int_{\mathbb R^n} f(y) G_s(x-y)\,dy &=
\int_0^\infty f(\rho) \rho^{n-1} \left( \int_{S^{n-1}} G_s(re_1 - \rho \omega) d\sigma_\omega\right)  d\rho,
\end{align*}
where $r=|x|$. In polar coordinates we can write
$$
|u(x)|\le\int_0^\infty|f(\rho)|\rho^{n-1}K(r,\rho)\,d\rho \qquad \text{ where }
K(r,\rho):=\int_{S^{n-1}}G_s(re_1-\rho\omega)\,d\sigma_\omega.
$$
We claim that
$$
K(r,\rho)\le C\,\max\{r,\rho\}^{1-n}\,|r-\rho|^{s-1}.
$$
Indeed, since $G_s(z)\le c_1|z|^{s-n}$, it suffices to bound $J(r,\rho)=\int_{S^{n-1}}|re_1-\rho\omega|^{s-n}d\sigma_\omega$.
 
\emph{Case $\rho\notin(r/2,2r)$.} Then $|re_1-\rho\omega|\ge|r-\rho|\ge\max\{r,\rho\}/2$, so $J\le C\max\{r,\rho\}^{s-n}$. Since $|r-\rho|\le\max\{r,\rho\}$ and $s-1<0$,
$$\max\{r,\rho\}^{s-n}=\max\{r,\rho\}^{1-n}\max\{r,\rho\}^{s-1}\le \max\{r,\rho\}^{1-n}|r-\rho|^{s-1}.$$
 
\emph{Case $\rho\in(r/2,2r)$.} Let $\theta\in[0,\pi]$ be the angle between $e_1$ and $\omega$. Then
$$|re_1-\rho\omega|^2=(r-\rho)^2+2r\rho(1-\cos\theta)\ \ge\ (r-\rho)^2+\tfrac{4}{\pi^2}\,r\rho\,\theta^2,$$
and $d\sigma\le C\theta^{n-2}d\theta$. Substituting $\theta=|r-\rho|\,t/\sqrt{r\rho}$,
\begin{equation}\label{ineq J}
J\le C\,(r\rho)^{-\frac{n-1}{2}}|r-\rho|^{s-1}\int_0^\infty(1+t^2)^{\frac{s-n}{2}}t^{n-2}\,dt\le C\,r^{1-n}|r-\rho|^{s-1}.
\end{equation}
The integral is finite because the integrand behaves like $t^{s-2}$ at infinity and $s<1$.

\medskip

Therefore, we obtain the pointwise inequality
\begin{align*}
|u(x)|&\le\int_0^\infty|f(\rho)|\rho^{n-1}K(r,\rho)\,d\rho=:I_1+I_2+I_3,
\end{align*}
where $I_1,I_2,I_3$ correspond to $\rho<r/2$, $r/2\le\rho\le2r$ and $\rho>2r$, respectively.

\medskip
\noindent\textbf{Bound for $I_1$}, (i.e., $\rho<r/2$). Here $\max\{r,\rho\}=r$ and $|r-\rho|\ge r/2$, so the kernel is bounded by $C\rho^{n-1}r^{1-n}r^{s-1}=C\rho^{n-1}r^{s-n}$, and
$$
I_1\le C\,r^{s-n}\!\int_{B_r}|f|\,dx\le 2C\,r^{s-n}\|f\|_{L^A(B_r)}\|\chi_{B_r}\|_{L^{\hat A}}\le C\,r^{s}A^{-1}(r^{-n})\|f\|_{L^A(B_r)},
$$
since $\|\chi_{B_r}\|_{L^{\hat A}}=1/\hat A^{-1}(1/|B_r|)\le|B_r|A^{-1}(1/|B_r|)\le C\,r^nA^{-1}(r^{-n})$ (by concavity of $A^{-1}$: $A^{-1}(\lambda t)\le\max\{1,\lambda\}A^{-1}(t)$).
 
 \medskip
 
\noindent\textbf{Bound for $I_3$} (i.e., $\rho>2r$). Here $\max\{r,\rho\}=\rho$ and $|r-\rho|\ge\rho/2$, so the kernel is $\le C\rho^{n-1}\rho^{1-n}\rho^{s-1}=C\rho^{s-1}$. Then, by H\"older's inequality in Orlicz spaces and \eqref{norma.ext},
$$
\begin{aligned}
I_3&\le C\!\int_{2r}^\infty|f(\rho)|\rho^{s-1}d\rho=C'\!\int_{|x|>2r}|f(x)||x|^{s-n}dx\\
&\le 2C'\|f\|_{L^A(\mathbb R^n\setminus B_{2r})}\big\||x|^{s-n}\big\|_{L^{\hat A}(\mathbb R^n\setminus B_r)}\\
&\le C\,r^sA^{-1}(r^{-n})\|f\|_{L^A(\mathbb R^n\setminus B_{2r})}.
\end{aligned}
$$

\noindent\textbf{Bound for $I_2$} (i.e., $r/2\le\rho\le2r$). By \eqref{ineq J} and \eqref{norma.mid},
\begin{equation}
\begin{split}
I_2& \leq C\int_{r/2}^{2r}|f(\rho)|\rho^{n-1}r^{1-n}|r-\rho|^{s-1}\,d\rho \\ & = Cr^{1-n}\int_{S_r}|f(y)||r-|y||^{s-1}\,dy \leq  2C\,r^{1-n}\|f\|_{L^A(S_r)}\|||r-|y||^{s-1}\|_{L^{\hat A}(S_r)}\\ & \leq Cr^{s}A^{-1}(r^{-n})\|f\|_{L^{A}(S_r)}.
\end{split}
\end{equation}
\medskip
Therefore, from the bounds for $I_1$, $I_2$ and $I_3$ the result follows for $s\in (0, 1)$.

Assume now that $s=1$. By \cite[Theorems 3.4 and 3.5]{OS}, $H^{1, A}(\mathbb R^n)=W^{1, A}(\mathbb R^n)$ with equivalent norms. By density, assume that $u$ is smooth and has compact support.   Write $u(x)=u(|x|)=u(r)$. Moreover, $u(R)=0$ for $R>0$ large enough.    Observe that $u(r)=-\int_r^R u'(t)\,dt$. Consider the measure $d\mu = t^{n-1}dt$. Then, H\"{o}lder's inequality  yields
\begin{align*}
|u(r)|&\leq \int_r^R |u'(t)|\,dt  = \int_r^R |u'(t)| t^{n-1} t^{1-n}\,dt=\int_r^R |u'(t)|  t^{1-n}\,d\mu\\
&\leq 2 \| u'\|_{L^A((0,R); d\mu)} \|t^{1-n}\|_{L^{\hat A}((r,R),d\mu)}.
\end{align*}
Using polar coordinates, we get that
$$
\omega_n\int_0^R A\left(\frac{|u'(t)|}{\lambda}\right)t^{n-1}\,dt = \int_0^R \int_{S^{n-1}} A\left(\frac{|u'(t)|}{\lambda}\right)t^{n-1}\, d\sigma dt = \int_{B_R} A\left(\frac{|\nabla u(x)|}{\lambda}\right)\,dx
$$
from where it is obtained that
$$
\int_0^R A\left(\frac{|u'(t)|}{\lambda}\right)t^{n-1}\,dt = 1 \iff \int_{B_R} \frac{1}{\omega_n}A\left(\frac{|\nabla u(x)|}{\lambda}\right)\,dx =1.
$$
This gives
$$
 \| u'\|_{L^A((0,R); d\mu)} = \| \nabla u\|_{L^{\frac{1}{\omega_n}A}(B_R)} \leq  c_n \| \nabla u\|_{L^{A}(B_R)}.
$$
Similarly,
$$\|t^{1-n}\|_{L^{\hat A}((r,R),d\mu)}= \||x|^{1-n}\|_{L^{\frac{1}{\omega_n}\hat A}(B_R\setminus B_r)}\leq c_n\||x|^{1-n}\|_{L^{\hat A}(B_R\setminus B_r).}$$

Hence, by \eqref{norma.ext}, 
$$|u(x)| \leq C\| \nabla u\|_{L^{A}(B_R)}\||x|^{1-n}\|_{L^{\hat A}(B_R\setminus B_r)} \leq C rA^{-1}(r^{-n})\|u\|_{H^{1, A}(\mathbb R^n)}$$
and the result follows for $s=1$ as well.
\end{proof}

 The following theorem is an application of the Ni-type inequality Theorem \ref{Ni Lemma} to compactness of the embedding between potential spaces and weighted Orlicz-Lebesgue spaces. 
 
  \begin{theorem} \label{teo.rad.compacto}
 Assume that $A$, $\hat A$, and $B$ are $N$-functions satisfying the $\Delta_2$-condition. Assume that $s\in (0,1]$ is such that $\frac{1}{p^{-}_A}<s<\frac{n}{p^+_A}$ and that
 \begin{equation}\label{assumption A and B near 0}
p_A^+<p_B^-.
\end{equation}

 Let $\alpha>-n$. Suppose further that
 \begin{equation}\label{assumpt alpha}
 (p_B^--p_A^+)\left(\frac{1-n}{p_{A}^+}\right)+\alpha<0.
 \end{equation}and $p_{B}^+< \frac{(\alpha+n)p_{A}^-}{n-sp_{A}^-}$.  Then, the embedding $H_{rad}^{s,A}(\mathbb R^n) \hookrightarrow L^B(\mathbb R^n, |x|^\alpha\,dx)$ is compact.
\end{theorem}

\begin{remark}Observe that, for appropriate $\alpha$, the $N$-function $B=A_{n/s}$ satisfies the above conditions  under the additional assumptions that $p^+<(p^-)^*$ and
\begin{equation}\label{add A}
\dfrac{n-sp^-}{n-sp^+}\dfrac{np^+}{p^-}<\dfrac{n-1}{n-sp^-}\dfrac{np^-}{p^+}+1.
\end{equation}
Indeed, we have
$$p_{B}^+=(p^{+})^{*}:=\frac{np^+}{n-sp^+} \quad \text{and }\quad p_{A}^-=p^{-}.$$Assumptions $p_{B}^+< \frac{(\alpha+n)p_{A}^-}{n-sp_{A}^-}$ and  \eqref{assumpt alpha} are possible if we can find $\alpha$ verifying
\begin{equation}
(p^{+})^{*}\left(\dfrac{n-sp^{-}}{p^{-}}\right)-n < \alpha < ((p^{-})^{*}-p^{+})\left(\dfrac{n-1}{p^{+}}\right).
\end{equation}As can be checked, this is true in view of \eqref{add A}.

\end{remark}

\begin{proof}
Let $u_k$ be a sequence in $H_{rad}^{s, A}(\mathbb R^n)$ such that $u_k\rightharpoonup 0$. Write for $0<\varepsilon<\delta\leq 1$,
\begin{equation}
\begin{split}
\int_{\mathbb R^n}B(|u_k(x)|)|x|^\alpha\,dx = \int_{|x|<\varepsilon}B(|u_k(x)|)|x|^\alpha\,dx+\int_{\varepsilon \leq |x|< 1/\varepsilon}B(|u_k(x)|)|x|^\alpha\,dx\\+\int_{|x|\geq 1/\varepsilon}B(|u_k(x)|)|x|^\alpha\,dx:=I_1+I_2+I_3.
\end{split}
\end{equation}

\noindent \textbf{Bound for $I_1$}. Observe that $sp^-_A<n$, then by Ni's Theorem \ref{Ni Lemma}, 
\begin{equation}\label{I_1}
\begin{split}
I_1 &\leq \int_{|x|<\varepsilon}B\left(C|x|^s A^{-1}(|x|^{-n})\|u_k\|_{H^{s, A}(\mathbb R^n)} \right)|x|^\alpha\,dx \leq C_n\int_0^\varepsilon B(r^sA^{-1}(r^{-n}))r^{\alpha+n-1}\,dr\\& \leq C_n\int_0^\varepsilon B(r^{s-\frac{n}{p_{A}^-}})r^{n+\alpha-1}\,dr \leq C_n\int_0^\varepsilon r^{p_{B}^+\left(s-\frac{n}{p_{A}^-} \right)+\alpha+n-1}\,dr.
\end{split}
\end{equation}
where we have used the inequality $A^{-1}(t)\leq t^{\frac{1}{p_{A}^-}}$ for $ t >1$. Since $p_{B}^+< \frac{(\alpha+n)p_{A}^-}{n-sp_{A}^-}$ and $\alpha>-n$, we get
$$p_{B}^+\left(s-\frac{n}{p_{A}^-} \right)+\alpha+n > \frac{(\alpha+n)p_{A}^-}{n-sp_{A}^-}\left(s-\frac{n}{p_{A}^-} \right)+\alpha+n =0.$$

Hence, from \eqref{I_1},
\begin{equation}\label{I}
I_1\leq C\varepsilon^{p_{B}^+\left(s-\frac{n}{p_{A}^-} \right)+\alpha+n}.
\end{equation}
 
\noindent \textbf{Bound for $I_2$}. First, observe that \eqref{assumption A and B near 0} implies
that  there is $t_0>0$ such that for all $0<t<t_0$, it follows $B(t)\leq A(t)$. Moreover, by Strauss Theorem \ref{Strauss Lemma} and the inequality 
\begin{equation} \label{desig.AA}
t\leq A^{-1}(t)\hat{A}^{-1}(t)
\end{equation}
we get for $\varepsilon \leq |x|<1/\varepsilon$ that 
$$
|u_k(x)|\leq C\dfrac{|x|^{1-n}}{\hat{A}^{-1}(||^{1-n})} \leq CA^{-1}(|x|^{1-n})\leq CA^{-1}(\varepsilon^{1-n})\leq C\varepsilon^{\frac{1-n}{p^{-}}}.
$$
So,  there is a constant $C(\varepsilon)>1$ such that for $t_0\leq t \leq C\varepsilon^{\frac{1-n}{p^{-}}}$, 
$$B(t)\leq C(\varepsilon)\tilde{B}(t)$$
where $\tilde{B}$ is any $N$-functions verifying the assumptions of Theorem  \ref{compact emb}. In this way, we get for $t \leq C\varepsilon^{\frac{1-n}{p^{-}}}$, that
$$B(t)\leq C(\varepsilon)(A(t)+\tilde{B}(t)).$$
Hence, by Theorem  \ref{compact emb}, and up to a subsequence that we do not re label, $u_k\to 0$ in $L^{\tilde{B}}(\mathbb R^n)$. On the other hand, since  $H^{s, A}_{rad}(\mathbb R^n)\subset W^{s', A}(\mathbb R^n)$ for $0<s'<s$ (\cite[Theorem 3.7]{OS}), we get that $u_k$ is bounded in $W^{s', A}(\mathbb R^n)$, $u_k\rightharpoonup 0$ in $W^{s', A}(\mathbb R^n)$ and, by \cite[Theorem 9.1]{ACPS}, $u_k\to 0$ in $L^A(B_{1/\varepsilon})$. As a result,
\begin{equation}\label{II}
I_2\leq C(\varepsilon)\varepsilon^{-|\alpha|}\int_{\varepsilon <|x|\leq 1/\varepsilon}\left(A(|u_k(x)|)+\tilde{B}(|u_k(x)|)\right)dx\to 0 \quad \text{as }k\to \infty.
\end{equation}

\noindent \textbf{Bound for $I_3$}. 
Observe that 
\begin{equation}\label{III}
I_3\leq C\int_{1/\varepsilon}^\infty B(|u_k(r)|)r^{n+\alpha-1}\,dr =   C\int_{1/\varepsilon}^\infty \dfrac{ B(|u_k(r)|)}{ A(|u_k(r)|)} A(|u_k(r)|)r^{n+\alpha-1}\,dr. 
\end{equation}Now,   condition \eqref{assumption A and B near 0} implies that $p_B^->p_A^+$. Hence, $B(t)/A(t)$ is increasing since its derivative satisfies
\begin{equation*}
\begin{split}
\dfrac{b(t)A(t)-B(t)a(t)}{A(t)^2}=\dfrac{tb(t)A(t)-B(t)ta(t)}{tA(t)^2}\geq \dfrac{(p_B^--p_A^+)A(t)B(t)}{tA(t)^2}\geq 0.
\end{split}
\end{equation*} 
Observe that using the Strauss' Lemma \ref{Strauss Lemma}  and \eqref{desig.AA} we get
$$
|u_k(r)|\leq A^{-1}(r^{1-n}) \|u\|_{H^{s,A}(\mathbb R^n)}, 
$$
which leads to 
\begin{align*}
\frac{ B(|u_k(r)|)}{ A(|u_k(r)|)} &\leq \dfrac{ B(A^{-1}(r^{1-n}) \|u_k\|_{H^{s,A}(\mathbb R^n)})}{ A(A^{-1}(r^{1-n}) \|u_k\|_{H^{s,A}(\mathbb R^n)})}\\
&\leq  \dfrac{\max\left\lbrace \|u_k\|_{H^{s,A}}^{p_B^+},\|u_k\|_{H^{s,A}(\mathbb R^n)}^{p_B^-} \right\rbrace}{\min\left\lbrace \|u_k\|_{H^{s,A}}^{p_A^+}, \|u_k\|_{H^{s,A}(\mathbb R^n)}^{p_A^-} \right\rbrace} \dfrac{B(A^{-1}(r^{1-n}))}{r^{1-n}} \leq C\dfrac{B(A^{-1}(r^{1-n}))}{r^{1-n}},
\end{align*}
where $C=C(\|u_k\|_{H^{s,A}(\mathbb R^n)}, p^\pm_A, p^\pm_B$) since $p_B^->p_A^+$. Hence,

\begin{equation}\label{IV}
\begin{split}
\int_{1/\varepsilon}^\infty \dfrac{ B(|u_k(r)|)}{ A(|u_k(r)|)} A(|u_k(r)|)r^{n+\alpha-1}\,dr &\leq C \int_{1/\varepsilon}^\infty \dfrac{B(A^{-1}(r^{1-n}))}{r^{1-n}}A(|u_k(r)|)r^{\alpha+n-1}\,dr \\ & \leq C\int_{1/\varepsilon}^\infty r^{p_B^-\left(\frac{1-n}{p_A^+}\right)-1+n+\alpha}A(|u_k(r)|)r^{n-1}\,dr \\& \leq C\varepsilon^{-\left[(p_B^--p_A^{+})\left(\frac{1-n}{p_A^+}\right)+\alpha\right]}\int_{\mathbb R^n}A(|u_k(x)|)\,dx,
\end{split}
\end{equation}where $(p_B^--p_A^{+})\left(\frac{1-n}{p_A^+}\right)+\alpha<0$ by assumption \eqref{assumpt alpha}. Hence, the strong convergence in $ L^B(\mathbb R^n, |x|^\alpha\,dx)$ follows from \eqref{I}, \eqref{II}, \eqref{III} and \eqref{IV}.
\end{proof}

\begin{theorem}\label{teo.sob.compacto}
Under the assumptions of Theorem \ref{teo.rad.compacto},  the embedding $W_{rad}^{s,A}(\mathbb R^n) \hookrightarrow L^B(\mathbb R^n, |x|^\alpha\,dx)$ is compact.
\end{theorem}

\begin{proof}
Let $0<s'<s$ be close enough to $s$ such that $\frac{1}{p_A^{-}}<s'$ and $p_B^{+}<\frac{(\alpha+n)p_A^{-}}{n-s'p_A^{-}}$. Then, according to Theorem \ref{teo.rad.compacto}, $H_{rad}^{s',A}(\mathbb R^n) \hookrightarrow L^B(\mathbb R^n, |x|^\alpha\,dx)$ compactly. Since by Theorem 3.9 in \cite{OS}, $W^{s, A}(\mathbb{R}^n)\hookrightarrow H^{s', A}(\mathbb{R}^n)$ continuously, we end the proof.
\end{proof}

\section{Application to a mixed local-nonlocal problem}\label{application}
The goal of this section is to extend existence results of \cite{OS2} for H\'enon-type equations in $\mathbb{R}^n$ to operators exhibiting both local and nonlocal behavior in the Orlicz setting. 

We assume that $A$ and $B$ are  $N$-functions satisfying the assumptions of Theorem \ref{teo.rad.compacto}. We consider 
\begin{equation} \label{eq}
\begin{cases}
\mathcal{L}_{s,A,\gamma}\, u + \dfrac{a(|u|)}{|u|}u =|x|^\alpha \dfrac{b(|u|)}{|u|}u  \quad \text{in }\mathbb{R}^n, \\
u >0,
\end{cases}
\end{equation}
where $u$ belongs to a suitable functional space of radial functions $\mathcal{W}^{s,A}(\mathbb{R}^n)$  (see \eqref{space}). Given $\gamma\in [0, 1]$, $s\in (0, 1]$ (depending on $\gamma$), the operator $\mathcal{L}_{s,A,\gamma}$ is given by 
\begin{equation*}
\mathcal{L}_{s,A,\gamma}\, u := \begin{cases}
(-\Delta_a) u, \quad & \text{if}\quad \gamma=1,\\
\gamma (-\Delta_a) u+ (1-\gamma)(-\Delta_a)^s u, \quad & \text{if}\quad \gamma\in (0,1), \, s\in (0,1),\\(-\Delta_a)^s u, \quad & \text{if}\quad \gamma=0,  \, s\in (0,1),
\end{cases}
\end{equation*}where we denote the fractional $a-$Laplacian of order $s\in(0,1]$ by
\begin{align*}
(-\Delta_a)^s u(x) :=
\begin{cases}
 2\,\text{p.v.}\displaystyle  \int_{\mathbb{R}^n} \dfrac{a(|D_su(x, y)|)}{|D_su(x, y)|}D_su(x, y) \dfrac{dy}{|x-y|^{n+s}}, & \text{ if } s\in (0,1)\\[15pt]
-\text{div}\,\left(\dfrac{a(|\nabla u|)}{|\nabla u|} \nabla u(x)\right), & \text{ if } s=1,
\end{cases}
\end{align*}where
$$D_su(x, y)=\dfrac{u(x)-u(y)}{|x-y|^s}.$$

Hence, the operator $\mathcal{L}_{s,A,\gamma}$ is a mixed local–nonlocal operator, reducing when $A(t)=t^p$ $t\geq 0$ to the classical 
$p-$Laplacian for  $\gamma=1$, and to the fractional $p-$Laplacian when $\gamma=0$. We will usually take $s\in (0, 1]$, identifying the case $\gamma=1$ with $s=1$.

Given $A$ an $N$-function and $s\in(0,1]$  we recall that 
$$
W^{s, A}(\mathbb{R}^n)= 
 \{ u \in L^A(\mathbb{R}^n) \colon [u]_{s,A}<\infty\}
$$
endowed with the norm  
$$
\|u\|_{W^{s, A}(\mathbb R^n)}:= \|u\|_{L^A(\mathbb R^n)} + [u]_{W^{s, A}(\mathbb R^n)},
$$
where 
$$
\|u\|_{L^A(\mathbb R^n)} :=\inf\left\{\lambda>0\colon \rho_A\left(\frac{u}{\lambda}\right)\leq 1 \right\},
$$
and ~the seminorm is given by
\begin{align*}
[u]_{W^{s, A}(\mathbb R^n)}:=
\begin{cases}\vspace{0.2cm}
\inf \left\lbrace \lambda>0:\displaystyle \rho_{s,A}\left(\frac{u}{\lambda}\right) \leq 1 \right\rbrace &\text{ if } s\in (0,1)\\ 
\inf \left\lbrace \lambda>0: \displaystyle   \rho_A\left(\frac{|\nabla u|}{\lambda}\right)\leq 1 \right\rbrace & \text{ if } s=1.
\end{cases} 
\end{align*}
where the modular functions are defined as follows
$$
\rho_{s,A}(u):=\int_{\mathbb R^n}\int_{\mathbb R^n} A(|D^s u|)\frac{dxdy}{|x-y|^n}, \quad \rho_A(u):=\int_{\mathbb R^n} A(|u|)\,dx.
$$

For Orlicz weighted spaces, we consider the norm
$$\|u\|_{L^{A,\alpha}(\mathbb R^n)} :=\inf\left\{\lambda>0\colon \rho_{A,\alpha}\left(\frac{u}{\lambda}\right)\leq 1 \right\}$$where
$$\rho_{A,\alpha}(u):= \int_{\mathbb R^n} A(|u|)|x|^\alpha\,dx.$$
We recall the following relation between norms and modular, provided that $A$ satisfies \eqref{eq.p}
\begin{equation} \label{relacion}
\|u\|_{L^A(\mathbb R^n)}^{p^+} \leq \rho_A(u) \leq \|u\|_{L^A(\mathbb R^n)}^{p^-} \quad \text{ if } \|u\|_{L^A(\mathbb R^n)} \leq 1, \qquad 
\|u\|_{L^A(\mathbb R^n)}^{p^-} \leq \rho_A(u) \leq \|u\|_A^{p^+} \quad \text{ if } \|u\|_{L^A(\mathbb R^n)} >  1.
\end{equation}
Moreover,  (see for instance \cite[Lemma 4.2]{FBS}), for $u\in W^{1,A}(\mathbb{R}^n)$,  and $s\in (0,1)$, it holds that
\begin{equation} \label{cota.local.nolocal}
\rho_{s,A}(u) \leq \dfrac{n \omega_n}{1-s} \rho_A(|\nabla u|)+ \dfrac{2Cn \omega_n}{s}\rho_A(u).
\end{equation}
Then,
\begin{equation} \label{cota.lnl}
[u]_{{W^{s,A}(\mathbb R^n)}}\leq C\|u\|_{W^{1,A}(\mathbb R^n)}.
\end{equation}

With the aforementioned definitions,  the natural framework for introducing weak solutions of \eqref{eq} is the space
\begin{align}\label{space}
\mathcal{W}^{s,A}(\mathbb{R}^n)=
\begin{cases}
W^{s,A}_{rad}(\mathbb{R}^n) & \text{ if } \gamma=0, s\in (0,1)\\
W^{1,A}_{rad}(\mathbb{R}^n)  & \text{ if } \gamma\in(0,1],
\end{cases}
\end{align}
endowed with 
$\|\cdot\|_{W^{s,A}(\mathbb R^n)}$ when $\gamma=0$ and $s\in (0, 1)$, and $\|\cdot\|_{W^{1,A}(\mathbb R^n)}$ when $\gamma \in (0, 1]$. 

We point out that since $W^{1, A}(\mathbb{R}^n)\subset W^{s, A}(\mathbb{R}^n)$ for all $s\in (0, 1)$ (see Lemma 4.2 in \cite{FBS}), the following representation formula holds: given $u\in \mathcal{W}^{s,A}(\mathbb{R}^n)$,
\begin{align}\label{representation}
\langle \mathcal{L}_{s,A,\gamma} \,u ,v \rangle &=\gamma \int_{\mathbb{R}^n}\dfrac{a(|\nabla u|)}{|\nabla u|}\nabla u \cdot \nabla v\,dx + (1-\gamma)  \int_{\mathbb{R}^n}\int_{\mathbb{R}^n} \dfrac{a(|D_s u(x, y)|)}{|D_s u(x, y)|}D_su(x, y)D_sv(x, y) \dfrac{dx\,dy}{|x-y|^n}
\end{align} 
for all $v\in \mathcal{W}^{s,A}(\mathbb{R}^n)$.  Therefore,  we say that $u\in \mathcal{W}^{s,A}(\mathbb{R}^n)$ is a \emph{weak solution} of \eqref{eq} if
$$
\langle \mathcal{L}_{s,A,\gamma}u,v \rangle  + \int_{\mathbb{R}^n} \dfrac{a(|u|)}{|u|}u v\,dx  = \int_{\mathbb{R}^n} |x|^\alpha \dfrac{b(|u|)}{|u|}u v\,dx \qquad \forall v\in \mathcal{W}^{s,A}(\mathbb{R}^n).
$$

The corresponding energy functional associated to \eqref{eq} is given by
$$
\mathcal{J}(u):= \gamma \rho_{A}(|\nabla u|)+ (1-\gamma)\rho_{s,A}(u) + \rho_A(u) - \rho_{B,\alpha}(u), \qquad \gamma \in [0,1].
$$
Indeed, $\mathcal{J}$ is a $C^1$ functional whose Fr\'echet derivatives is given by
$$
\langle \mathcal{J}'(u),v\rangle  = \langle \mathcal{L}_{s,A,\gamma} \,u ,v \rangle  +\int_{\mathbb R^n} \frac{a(|u|)}{|u|}uv\,dx - \int_{\mathbb R^n} \frac{b(|u|)}{|u|}uv |x|^\alpha\,dx \qquad \forall u,v \in \mathcal{W}^{s,A}(\mathbb{R}^n). 
$$

Our main application concerns the existence of radial solutions of problem \eqref{eq}.

\begin{theorem}\label{existence theorem}
Assume that $A$, $\hat A$, and $B$ are $N$-functions satisfying the $\Delta_2$-condition and suppose that $A(\sqrt{t})$ is convex. Assume that $\frac{1}{p^{-}_A}<s<\frac{n}{p^+_A}$ and that $p_A^+<p_B^-$.  Let $\alpha>-n$ such that
 $$
 (p_B^--p_A^+)\left(\frac{1-n}{p_{A}^+}\right)+\alpha<0, \qquad \text{and} \qquad p_{B}^+< \frac{(\alpha+n)p_{A}^-}{n-sp_{A}^-}.
 $$
Then, problem \eqref{eq} admits a radial weak solution. 
\end{theorem}
\begin{proof}
We check the conditions of the Brezis-Nirenberg Lemma  \ref{lema.BN} with  $E=\mathcal{W}^{s,A}(\mathbb{R}^n)$.

\textbf{(i)}\, Case $\gamma=0$ and $s\in (0,1)$. Here the underlying space is $W^{s,A}_{rad}(\mathbb R^n)$. Let $u$  be such that $\|u\|_{W^{s,A}(\mathbb R^n)}\leq 1$.  By Theorem \ref{teo.rad.compacto} we can assume that $\|u\|_{W^{s,A}(\mathbb R^n)}\leq 1$ is small enough such that $\|u\|_{B,\alpha}\leq 1$. Then, using \eqref{relacion} and Theorem \ref{teo.rad.compacto}
\begin{align*}
\mathcal{J}(u)&= \rho_{s,A}(u)+ \rho_A(u) - \rho_{B,\alpha}(u)\geq \|u\|_{W^{s,A}(\mathbb R^n)}^{p^+_A} - \|u\|_{L^{B,\alpha}(\mathbb R^n)}^{p^-_B}
\geq \|u\|_{W^{s,A}(\mathbb R^n)}^{p^+_A} - C\|u\|_{W^{s,A}(\mathbb R^n)}^{p^-_B}
\end{align*}
Since $p^+_A < p^-_B$, we have that $\mathcal{J}(u)>0$ for $\|u\|_{W^{s,A}(\mathbb R^n)}=r$ small enough.

Case $\gamma=1$ ($s=1$). Let $u\in W^{1,A}_{rad}(\mathbb R^n)$ be such that $\|u\|_{W^{1,A}(\mathbb R^n)}\leq 1$. Similarly, we can assume that $\|u\|_{B,\alpha}\leq 1$. From \eqref{relacion} and Theorem \ref{teo.rad.compacto}
\begin{align*}
\mathcal{J}(u)&= \rho_{A}(|\nabla u|)+ \rho_A(u) - \rho_{B,\alpha}(u)\geq \|u\|_{W^{1,A}(\mathbb R^n)}^{p^+_A} - \|u\|_{L^{B,\alpha}(\mathbb R^n)}^{p^-_B}
\geq \|u\|_{W^{1,A}(\mathbb R^n)}^{p^+_A} - C\|u\|_{W^{1,A}(\mathbb R^n)}^{p^-_B}
\end{align*}
Since $p^+_A < p^-_B$, we have that $\mathcal{J}(u)>0$ for $\|u\|_{W^{1,A}(\mathbb R^n)}=r$ small enough.

Case $\gamma\in (0,1)$ and $s\in (0,1)$. Let $u\in W^{1,A}_{rad}(\mathbb R^n)$. By Theorem \ref{teo.rad.compacto} we can assume that $\|u\|_{W^{1,A}(\mathbb R^n)}\leq 1$ is small enough such that $\|u\|_{L^{B,\alpha}(\mathbb R^n)}\leq 1$.  Then, using \eqref{relacion} and Theorem \ref{teo.rad.compacto}
\begin{align*}
\mathcal{J}(u)&= \gamma \rho_{A}(|\nabla u|)+ (1-\gamma)\rho_{s,A}(u) + \rho_A(u) - \rho_{B,\alpha}(u)\\
&\geq \gamma \|u\|_{W^{1,A}(\mathbb R^n)}^{p^+_A} - \|u\|_{L^{B,\alpha}(\mathbb R^n)}^{p^-_B}
\geq \gamma \|u\|_{W^{1,A}(\mathbb R^n)}^{p^+_A} - C\|u\|_{W^{1,A}(\mathbb R^n)}^{p^-_B}
\end{align*}
Since $p^+_A < p^-_B$, we have that $\mathcal{J}(u)>0$ for $\|u\|_{W^{1,A}(\mathbb R^n)}=r$ small enough.

\medskip

\noindent \textbf{(ii)}\,  Let $\gamma \in [0,1]$ and $s\in (0,1]$ (where we identify $\gamma=1$ with $s=1$). Let $w>0$ fixed, and define $v:=tw$, $t\geq 1$. Then, using \eqref{relacion}
\begin{align*}
\mathcal{J}(v)&= \gamma \rho_{A}(t|\nabla w|)+ (1-\gamma)\rho_{s,A}(tw) + \rho_A(tw) - \rho_{B,\alpha}(tw)\\
&\leq 
 t^{p^+_A}\gamma \rho_{A}(|\nabla w|)+ t^{p^+_A}(1-\gamma)\rho_{s,A}(w) + t^{p^+_A}\rho_A(w) -t^{p^-_B} \rho_{B,\alpha}(w)\\
 &\leq C_1 t^{p^+_A} - C_2 t^{p^-_B}
\end{align*}
for $C_1$ and $C_2$ positive constants depending on $w$, $\gamma$ and $A$. Since $p^+_A<p^-_B$, there exists $t$ such that $\mathcal{J}(v)<0$ and $\|v\|>r$.

\medskip

\noindent \textbf{(iii)}\,  Let $\gamma\in (0,1)$ and $s\in (0,1)$. Let  $\{u_k\}_{k\in\mathbb{N}}\subset W^{1,A}_{rad}(\mathbb{R}^n)$  be such that $\mathcal{J}'(u_k)\to 0$ and $\mathcal{J}(u_k)\to c$. Since $\mathcal{J}'(u_k)\to 0$, using \eqref{cota.lnl} we get that
\begin{align*}
|\langle \mathcal{J}'(u_k),u_k\rangle|  &= \left|\langle \mathcal{L}_{s,A,\gamma} \,u_k ,u_k \rangle  + \int_{\mathbb R^n} a(|u_k|)|u_k|\,dx - \int_{\mathbb R^n} b(|u_k|)|u_k| |x|^\alpha\,dx \right|\\
&\leq
\| \mathcal{J}'(u_k) \|_{ W^{1,A}_{rad}(\mathbb R^n)'} \|u_k\|_{1,A}\\
&\leq
 \|u_k\|_{W^{1,A}(\mathbb R^n)}
\end{align*}
for $k$ large enough. In particular, if we define
$$
E_k:= \gamma  \rho_A(|\nabla u_k|) + (1-\gamma) \rho_{s,A}(u_k) +  \rho_A(u_k),
$$
for $k$ large enough we have that
\begin{equation} \label{eq.t.1}
p^-_B   \rho_{B,\alpha}(u_k) \leq  \|u_k\|_{W^{1,A}(\mathbb R^n)} + p^+_A E_k.
\end{equation}
Since $|\mathcal{J}(u_k)|\leq C$, we have that
\begin{equation} \label{eq.t.2}
|\mathcal{J}(u_k)|= \left| E_k  - \rho_{B,\alpha}(u_k)\right| \leq C.
\end{equation}
From \eqref{eq.t.1} and \eqref{eq.t.2} we get that
\begin{align*}
E_k &  \leq C + \rho_{B,\alpha}(u_k) \leq C + 
\frac{1}{p^-_B}   \|u_k\|_{W^{1,A}(\mathbb R^n)} + \frac{p^+_A}{p^-_B}E_k
\end{align*}
from where, since $p^+_A<p^-_B$, 
\begin{align*}
\left(1-\frac{p^+_A}{p^-_B}\right)E_k \leq  C + 
\frac{1}{p^-_B}   \|u_k\|_{W^{1,A}(\mathbb R^n)}
\end{align*}
This gives that for some $c>0$ depending on $p^+_A$ and $p^-_B$
$$
E_k \leq  C_1(1 + 
 \|u_k\|_{W^{1,A}(\mathbb R^n)})
$$
if $\|u_k\|_{W^{1,A}(\mathbb R^n)}>1$,  using \eqref{relacion}
$$
\|u_k\|_{W^{1,A}(\mathbb R^n)}^{p^-_A} \leq 
C_1(1+\|u_k\|_{W^{1,A}(\mathbb R^n)}) \implies \|u_k\|_{W^{1,A}(\mathbb R^n)}\leq (2C_1)^\frac{1}{p^-_A -1}.
$$
Therefore, it follows that $\{u_k\}_{k\in\mathbb{N}}$ is bounded in $W^{1,A}_{rad}(\mathbb{R}^n)$.

Then, up to a subsequence, $u_k\rightharpoonup u$ weakly in $W^{1,A}_{rad}(\mathbb{R}^n)$. From the compact embedding given in Theorem \ref{teo.rad.compacto}, we get that
$u_k \to u$ strongly in $L^B(\mathbb{R}^n, |x|^\alpha dx)$, which implies, due to convexity that
$$
\int_{\mathbb{R}^n} |x|^\alpha b(u_k) (u_k-u) \,dx \to 0
$$
as $k\to\infty$. 
Therefore
\begin{align*}
\Big|\langle \mathcal{L}_{s,A,\gamma}\, u_k,u_k-u \rangle +\int_{\mathbb{R}^n} \frac{a(|u_k|)}{|u_k|}&u_k (u_k-u)\,dx  \Big| = \left|\langle \mathcal{J}'(u_k), u_k-u \rangle + \int_{\mathbb{R}^n} \frac{b(|u_k|)}{|u_k|}u_k (u_k-u)|x|^\alpha \,dx  \right|\\
&\leq
\|\mathcal{J}'(u_k)\|_{W^{1,A}(\mathbb{R}^n)'} \|u_k-u\|_{W^{1,A}(\mathbb R^n)} +
 \int_{\mathbb{R}^n} \frac{b(|u_k|)}{|u_k|}u_k (u_k-u)|x|^\alpha \,dx \to 0
 \end{align*}
 as $k\to\infty$, where $W^{1,A}(\mathbb{R}^n)'$ denotes the dual space of $W^{1,A}(\mathbb{R}^n)$.

Therefore,   the operator $\mathcal{\tilde L}$ defined  as
$$
\langle  \mathcal{\tilde L} u,v\rangle := \langle \mathcal{L}_{s,A,\gamma}\, u,v\rangle + \int_{\mathbb{R}^n} \frac{a(|u|)}{|u|}uv\,dx 
$$
fulfills that 
$$
\langle  \mathcal{\tilde L} u_k,u_k-u \rangle \to 0 \text{ as } k\to \infty.
$$
Hence, by Lemma \ref{rem.S.prop} below, we get that $u_k\to u$ strongly in $W^{1,A}(\mathbb{R}^n)$.

When $\gamma=0$ and $s\in (0,1)$, and when $\gamma=1$ (identifying this case with $s=1$), the reasoning is analogous.

\medskip

\noindent \textbf{(iv)}\,   Let $\gamma \in (0,1)$ and $s\in (0,1)$. To obtain a positive solution, we choose the mapping $P$ defined as $P(u)=|u|$. Trivially, $P(0)=0$ and $P(v)=v$, where $v$ is given (ii). Moreover, it is direct that $\| \nabla|u|\|_{L^A(\mathbb R^n)} = \|\nabla u\|_{L^A(\mathbb R^n)}$, $\||u|\|_{L^A(\mathbb R^n)}= \|u\|_{L^A(\mathbb R^n)}$  and $\||u|\|_{L^{A,\alpha}(\mathbb R^n)} = \|u\|_{L^{A,\alpha}(\mathbb R^n)}$. Furthermore, since $||u(x)|-|u(y)|| \leq |u(x)-u(y)|$, and $t\mapsto A(t)$ is monotone increasing for $t\geq 0$, $A(|D^s |u||) \leq A(|D^s u|)$. This gives that $\rho_{s,A}(|u|) \leq \rho_{s,A}(u)$ and hence $[|u|]_{W^{s,A}(\mathbb R^n)} \leq [u]_{W^{s,A}(\mathbb R^n)}$. All these relations give that for any $u\in W^{1,A}(\mathbb{R}^n)$,
$$
\mathcal{J}(P(u)) \leq \mathcal{J}(u).
$$
An analogous (and simpler) analysis can be made in the cases $\gamma=0$ and $s\in (0,1)$, and $\gamma=1$ (identified with $s=1$).

Therefore, the conclusion follow from  Lemma \ref{lema.BN}. This concludes the proof. 
\end{proof}

\begin{lemma}\label{rem.S.prop}Suppose that $A$ is an $N$-function satisfying the $\Delta_2$-condition and  suppose that $A(\sqrt{t})$ is convex. Let $\{u_k\}_{k\in\mathbb{N}} \in \mathcal{W}^{s, A}(\mathbb{R}^n)$ be  a sequence such that $u_k\rightharpoonup u$ and 
$$
\langle  \mathcal{\tilde L} u_k,u_k-u \rangle \to 0 \text{ as } k\to \infty, \qquad \langle  \mathcal{\tilde L} u,v\rangle := \langle \mathcal{L}_{s,A,\gamma}\, u,v\rangle + \int_{\mathbb{R}^n} \dfrac{a(|u|)}{|u|} uv\,dx .
$$where  $v$ is defined in Theorem \ref{existence theorem}. Then $u_k\to u$ strongly in $\mathcal{W}^{s, A}(\mathbb{R}^n)$.

\end{lemma}

\begin{proof}Observe that the weak convergence implies that
\begin{equation}\label{weak S}
\langle  \mathcal{\tilde L} u_k-\mathcal{\tilde L}u,u_k-u \rangle \to 0 \text{ as } k\to \infty.
\end{equation}
We will appeal to the following well-known inequality for $\xi, \eta \in \mathbb{R}^n$ (see \cite[Lemma 8]{O}):
\begin{equation}\label{ineq conv}
\left(\dfrac{a(|\xi|)}{|\xi|}\xi - \dfrac{a(|\eta|)}{|\eta|}\eta\right)(\xi-\eta)\geq CA(|\xi-\eta|).
\end{equation}
First, assume $\gamma=0$ and hence $\mathcal{W}^{s, A}(\mathbb{R}^n)=W^{s, A}_{rad}(\mathbb{R}^n).$ Then
\begin{equation*}
\begin{split}
& \langle  \mathcal{\tilde L} u_k-\mathcal{\tilde L}u,u_k-u \rangle \\& = \int_{\mathbb{R}^n}\int_{\mathbb{R}^n}\left(\dfrac{a(|D_su_k(x, y)|)}{|D_su_k(x, y)|}D_su_k(x, y)-\dfrac{a(|D_s u(x, y)|)}{|D_su(x, y)|}D_su(x, y) \right)\left(D_su_k(x, y)-D_su(x,y)\right)\dfrac{dx\,dy}{|x-y|^n}\\ & \quad +\int_{\mathbb{R}^n}\left(\dfrac{a(|u_k(x)|)}{|u_k(x)|}u_k(x)- \dfrac{a(|u(x)|)}{|u(x)|}u(x)\right)(u_k(x)-u(x))\,dx\\ &  
\geq C
\left(\rho_{s, A}(u_k-u)+\rho_A(u_k-u)\right).
\end{split}
\end{equation*}Hence, \eqref{weak S} gives the strong convergence in $W^{s, A}(\mathbb{R}^n)$. If now $\gamma\in (0, 1]$, then  $\mathcal{W}^{s, A}(\mathbb{R}^n)=W^{1, A}_{rad}(\mathbb{R}^n).$ Then, as before
\begin{equation}\label{similar i}
\begin{split}
 \langle  \mathcal{\tilde L} u_k-\mathcal{\tilde L}u,u_k-u \rangle \geq  C\left(\gamma\rho_{A}(\nabla u_k-\nabla u)+\rho_A(u_k-u)\right),
\end{split}
\end{equation}and the conclusion holds.
\end{proof}

\section*{Appendix}
The following variant of the Mountain-Pass Theorem will be of use for our arguments. See \cite[Theorem 10]{BN}. 
\begin{lemma}[Brezis--Nirenberg] \label{lema.BN}
Let $E$ be a Banach space and $\Phi \in C^{1}(E,\mathbb{R})$ satisfy $\Phi(0)=0$ and the following conditions:
\begin{itemize}
    \item[(i)] there exists $r>0$ such that $\displaystyle \inf_{\|u\|=r}\Phi(u)>0$;
    \item[(ii)] there exists $v$ with $\|v\|>r$ such that $\Phi(v)<0$;
    \item[(iii)] the Palais--Smale condition: if $\Phi'(u_k)\to 0$ and $\Phi(u_k)\to c$,  
    then $\{u_k\}_{k\in \mathbb{N}}$ has a subsequence converging to $u\in E$;
    \item[(iv)] there exists a continuous mapping $P:E\to E$ such that  
    $P(0)=0$, $P(v)=v$, and $\Phi(P(u))\le \Phi(u)$ for every $u\in E$.
\end{itemize}
Then there exists a critical point $u^{*}\in \overline{P(E)}$ such that
\[
\Phi(u^{*}) \ge \inf_{\|u\|=r} \Phi(u).
\]
\end{lemma}

\end{document}